\documentclass{amsart}

 \theoremstyle{definition}
 
 \theoremstyle{remark}
 
 \numberwithin{equation}{subsection}
\begin{document}

\title[2-LOCAL DERIVATIONS ON AW$^*$-ALGEBRAS]
 {2-LOCAL DERIVATIONS ON AW$^*$-ALGEBRAS OF TYPE I}

\author{Shavkat Ayupov}

\address{Institute of Mathematics, National University of Uzbekistan,
Tashkent, Uzbekistan and the Abdus Salam International Centre for
Theoretical Physics (ICTP) Trieste, Italy}

\email{sh$_-$ayupov@mail.ru}

\author{Farhodjon Arzikulov}

\address{Andizhan State University, Andizhan, Uzbekistan}

\email{arzikulovfn@rambler.ru}

\begin{abstract}
It is proved that every 2-local derivation on an
AW$^*$-algebra of type I is a derivation. Also an analog of
Gleason theorem for signed measures on projections of homogenous
AW$^*$-algebras except the cases of an AW$^*$-algebra of type
I$_2$ and a factor of type I$_m$, $2<m<\infty$ is proved.
\end{abstract}

\maketitle



\section*{Introduction}

The present paper is devoted to 2-local derivations on
AW$^*$-algebras. Recall that a 2-local derivation is defined as
follows: given an algebra $A$, a map $\Delta : A \to A$ (not
linear in general) is called a 2-local derivation if for every
$x$, $y\in A$, there exists a derivation $D_{x,y} : A\to A$ such
that $\Delta(x)=D_{x,y}(x)$ and $\Delta(y)=D_{x,y}(y)$.

In 1997, P. \v{S}emrl \cite{S} introduced the notion of 2-local
derivations and described 2-local derivations on the algebra
$B(H)$ of all bounded linear operators on the infinite-dimensional
separable Hilbert space H. A similar description for the
finite-dimensional case appeared later in \cite{KK}. In the paper
\cite{LW} 2-local derivations have been described on matrix
algebras over finite-dimensional division rings.

In \cite{AK} the authors suggested a new technique and have
generalized the above mentioned results of \cite{S} and \cite{KK}
for arbitrary Hilbert spaces. Namely they considered 2-local
derivations on the algebra $B(H)$ of all linear bounded operators
on an arbitrary (no separability is assumed) Hilbert space $H$ and
proved that every 2-local derivation on $B(H)$ is a derivation.

In \cite{AA} we also suggested another technique and generalized
the above mentioned results of \cite{S}, \cite{KK} and \cite{AK}
for arbitrary von Neumann algebras of type I and proved that every
2-local derivation on these algebras is a derivation. In
\cite{AKNA} (Theorem 3.4) a similar result was proved for finite
von Neumann algebras. In \cite{AA2} we extended all above results
and give a short proof of the theorem for arbitrary semi-finite
von Neumann algebras. Finally, in \cite{AK2} there was given a
proof of the problem for von Neumann algebras.

In the present paper we prove that every 2-local derivation on an
AW$^*$-algebra of type I is a derivation (theorem 3.1). Also we
prove an analog of Gleason theorem for signed measures on
projection of homogenous AW$^*$-algebras except the cases of an
AW$^*$-algebra of type I$_2$ and a factor of type I$_m$,
$2<m<\infty$ (theorem 1.2). Our proof is essentially based on this
analog of Gleason theorem for signed measures on projection of
homogenous AW$^*$-algebras.

\section*{Preliminaries}

Let $\mathcal{A}$ be an AW$^*$-algebra.

{\it Definition.} A linear map $D : \mathcal{A}\to \mathcal{A}$ is
called a derivation, if $D(xy)=D(x)y+xD(y)$ for any two elements
$x$, $y\in \mathcal{A}$.

A map $\Delta : \mathcal{A}\to \mathcal{A}$ is called a 2-local
derivation, if for any two elements $x$, $y\in \mathcal{A}$ there
exists a derivation $D_{x,y}:\mathcal{A}\to \mathcal{A}$ such that
$\Delta (x)=D_{x,y}(x)$, $\Delta (y)=D_{x,y}(y)$.

It is known that any derivation $D$ on a AW$^*$-algebra
$\mathcal{A}$ is an inner derivation \cite{Ol}, that is there
exists an element $a\in \mathcal{A}$ such that
$$
D(x)=ax-xa, x\in \mathcal{A}.
$$
Therefore for an AW$^*$-algebra $\mathcal{A}$ the above definition
is equivalent to the following one: A map $\Delta : \mathcal{A}\to
\mathcal{A}$ is called a 2-local derivation, if for any two
elements $x$, $y\in \mathcal{A}$ there exists an element $a\in
\mathcal{A}$ such that $\Delta (x)=ax-xa$, $\Delta (y)=ay-ya$.

Let $\mathcal{A}$ be an AW$^*$-algebra, $\Delta :\mathcal{A}\to
\mathcal{A}$ be a 2-local derivation. Then from the definition it
easily follows that $\Delta$ is homogenous. At the same time,
$$
\Delta(x^2)=\Delta(x)x+x\Delta(x)
$$
for each $x\in \mathcal{A}$.

In \cite{Bre} it is proved that any Jordan derivation on a
semi-prime algebra is a derivation. Since $\mathcal{A}$ is
semi-prime, the map  $\Delta$ is a derivation if it is additive.
Therefore, in the case of AW$^*$-algebra to prove that the 2-local
derivation $\Delta :\mathcal{A}\to \mathcal{A}$ is a derivation it
is sufficient to prove that $\Delta :\mathcal{A}\to \mathcal{A}$
is additive.

\medskip

\section{Gleason theorem and its application}

\medskip

{\it Definition.} Let $\mathcal{A}$ be an AW$^*$-algebra. The
lattice of all projections of $\mathcal{A}$ we denote by
$P(\mathcal{A})$. Recall that a map $\mu : P(\mathcal{A})\to {\Bbb
C}$ is called a signed measure (or charge) if
$\mu(e_1+e_2)=\mu(e_1)+\mu(e_2)$ for arbitrary mutually orthogonal
projections $e_1$, $e_2$ in $\mathcal{A}$.

A signed measure $\mu$ is said to be bounded if $\sup\{\vert
\mu(e)\vert : e\in P(\mathcal{A})\}$ is finite.

Recall that a map $\mu : P(\mathcal{A})\to {\Bbb R}_+$ is called a
finitely additive measure if $\mu(e_1+e_2)=\mu(e_1)+\mu(e_2)$ for
arbitrary mutually orthogonal projections $e_1$, $e_2$ in
$\mathcal{A}$, if additionally $\mu(1)=1$ then $\mu$ is called a
probability measure.

Let $A$ be a C*-algebra. By a singly generated C$^*$-subalgebra of
$A$ we mean a norm-closed $*$-subalgebra $A(x)$ generated by a
single self-adjoint element $x\in A$ (and the identity $1$ if $A$
has identity).

{\it Definition.} A positive quasi-linear functional is a function
$\rho: A\to {\Bbb C}$ such that

(i) $\rho\vert_{A(x)}$ is a positive linear functional for each
$x\in A_{sa}$.

(ii) $\rho(a)=\rho(a_1)+i\rho(a_2)$, when $a=a_1+ia_2$ is the
canonical decomposition of $a$ in self-adjoint parts $a_1$, $a_2$.

If in addition (iii) $\rho(1)=1$, then we say that $\rho$ is a
quasi-state on $A$.

Since an AW$^*$-factor of type {\rm I}$_n$ is a von Neumnn
algebra, we have

{\bf  Theorem 1.1.}  {\it Let $A$ be an AW$^*$-factor of type {\rm
I}$_n$, $n\neq 2$. Then every bounded signed measure on $P(A)$ can
be extended uniquely to a linear functional on $A$.}

Let $\mathcal{A}$ be an AW$^*$-algebra of type {\rm I}$_n$ that is
not a factor for $2<n<\infty$ and $n\neq 2$, where $n$ is a
cardinal number. Let $\{e_i\}$ be a maximal family of pairwise
orthogonal abelian projections in $\mathcal{A}$ with central
support 1 such that $\sup_i e_i=1$. Let $\{e_{ij}\}$ be the system
of matrix units with respect to $\{e_i\}$. Let $Z(\mathcal{A})$ be
the center of $\mathcal{A}$, $X$ be a compact such that $C(X)\cong
Z(\mathcal{A})$. Let $m$ be a natural number such that $m\leq n$,
$\{e_j\}_{j=1}^m$ be a subset of $\{e_i\}$ and
$e=\sum_{j=1}^me_j$. It is known that $e\mathcal{A}e\cong
C(X)\otimes M_m({\bf C})$. In \cite{Kus} it is proved that
$\mathcal{A}$ is isomorphic to the AW$^*$-algebra of maps from the
extremely disconnected compact $X$ to the von Neumann algebra
$B(H)$ of all bounded linear operators on a Hilbert space $H$
satisfying to certain conditions. Following by \cite{Kus} this
AW$^*$-algebra we denote by $SC_{\#}(X,B(H))$. Further we will use
these designations. The algebra $C(X,B(H))$ of all continuous
operator-valued functions on $X$ is a C$^*$-algebra.

{\bf  Lemma 1.2.}  {\it Every signed measure $\mu$ on
$P(SC_{\#}(X,B(H)))$ can be uniquely extended to a linear
functional on $C(X,B(H))$, i.e. there exists a linear functional
$\rho$ on $C(X,B(H))$ such that
$$
\rho\vert_{P(SC_{\#}(X,B(H)))\cap
C(X,B(H))}=\mu\vert_{P(SC_{\#}(X,B(H)))\cap C(X,B(H))}.
$$
}

{\bf Proof.} Let $\{E_i\}$ be a maximal family of pairwise
orthogonal Abelian projections in $B(H)$ and let $e_i:X\to E_i$ be
the operator-valued function $e_i(x)=E_i$, $x\in X$ for every $i$.
Then $\{e_i\}$ be a maximal family of pairwise orthogonal Abelian
projections in $SC_{\#}(X,B(H))$ with central support 1 such that
$\sup_i e_i=1$. Let $\{e_{ij}\}$ be the system of matrix units
with respect to $\{e_i\}$. The subalgebra
$$
\mathcal{C}=\{a\in A: \text{for all}\,\,\, i\,\,\,
\text{and}\,\,\, j\,\,\, e_iae_j=\lambda e_{ij}, \lambda\in {\Bbb
C}\}
$$
is an AW$^*$-subfactor of type {\rm I}$_n$.

It is clear that $\mathcal{C}\subset C(X,B(H)$.

Let $\mu$ be a signed measure on the set $\mathcal{P}$ of all
projections of $SC_{\#}(X,B(H))$. First, we will prove that the
signed measure $\mu$ is uniquely extended to a linear functional
on the vector space $C(X,B(H)_{sa})$ of all self-adjoint elements
of $C(X,B(H))$.

Let $\mu_{Re}(p)=Re(\mu(p))$, $\mu_{Im}(p)=Im(\mu(p))$ for every
$p\in P(\mathcal{A})$. Then
$$
\mu =\mu_{Re}+i\mu_{Im}.
$$
The maps $\mu_{Re}$, $\mu_{Im}$ are real number-valued signed
measures. Without loss of generality we assume that $\mu$ is a
real number-valued signed measures.

By theorem 1.1 the signed measure $\mu$ is uniquely extended to a
linear functional on any subalgebra of the form $z\mathcal{C}$,
where $z$  is a central projection.  This linear functional on
$z\mathcal{C}$ we denote by $\phi_z$. Clearly
$\phi_z\vert_{P(zC)}=\mu$.

At the same time, straightforward arguments show that $\mu$ has a
unique extension to a function $\delta: SC_{\#}(X,B(H))\to {\Bbb
C}$, where $\delta$ is linear and bounded on each abelian
$*$-subalgebra of $SC_{\#}(X,B(H))$ and where $\delta(x + iy)
=\delta(x) + i\delta(y)$ whenever $x$ and $y$ are self-adjoint
(see \cite{Aar}). Hence, for every central projection $z$ we have
$\phi_z\vert_{z\mathcal{C}}=\delta\vert_{z\mathcal{C}}$. Therefore
we may take
$$
\phi(za)=\phi_z(za), z\in Z(SC_{\#}(X,B(H)))\cap
P(SC_{\#}(X,B(H))), a\in \mathcal{C}.
$$
Clearly $\phi$ is a mapping.

Take the set $\Re$ of all elements of the form $\sum_{i=1}^{m}
x_ia_i$, where $x_1,\dots,x_{m}$ are orthogonal central
projections of sum 1 and the elements $a_1,\dots,a_{m}$ belong to
$\mathcal{C}_{sa}$. It is clear that $\Re$ is a normed space and
belongs to $C(X,B(H)_{sa})$.

For every element $\sum_{i=1}^{m} x_ia_i$ from $\Re$  we suppose
$$
\phi(\sum_{i=1}^{m} x_ia_i)= \sum_{i=1}^{m} \phi (x_ia_i).
$$

Then $\phi$  is a linear mapping on $\Re$. Indeed, for arbitrary
elements $\sum_{i=1}^{m_1} x_ia_i$ and $\sum_{j=1}^{m_2} y_jb_j$
of $\Re$ we have
$$
\phi (\sum_{i=1}^{m_1} a_iz_i+\sum_{j=1}^{m_2} y_jb_j)= \phi
[\sum_{i=1}^{m_1} \sum_{j=1}^{m_2} x_iy_j(a_i+b_j)]=
$$

$$
      \sum_{i=1}^{m_1} \sum_{j=1}^{m_2} \phi (x_iy_ja_i+
      x_iy_jb_j)=\sum_{i=1}^{m_1}\sum_{j=1}^{m_2}\phi (x_iy_ja_i)+
$$
$$
      \sum_{i=1}^{m_1} \sum_{j=1}^{m_2}\phi (x_iy_jb_j)=
      \phi(\sum_{i=1}^{m_1} x_ia_i)+\phi(\sum_{i=1}^{m_2} y_jb_j).
$$
Thus $\phi $ is a linear mapping on $\Re$.

Now we prove that, if a consequence $(a_n)\subset \Re$ is a
fundamental consequence, then the numerical consequence
$(\phi(a_n))$ is also fundamental.

By \cite[proposition 1.17.1]{S} there exist positive functionals
$\phi_z^+$, $\phi_z^-$ on $z\mathcal{C}$ for every central
projection $z$ such that
$$
\phi=\phi_z^+-\phi_z^-.
$$
At the same time, the measure
$$
\mu\vert_{Z(SC_{\#}(X,B(H)))\cap P(SC_{\#}(X,B(H)))}
$$
has a unique extension to a linear functional
$\eta :Z(SC_{\#}(X,B(H)))\to {\Bbb C}$, and by \cite[proposition
1.17.1]{S} there exist positive functionals $\eta^+$, $\eta^-$ on
$Z(SC_{\#}(X,B(H)))$ such that
$$
\eta=\eta^+-\eta^-.
$$
If $\lambda_z=\phi_z^+(z)-\eta^+(z)$ for some number $\lambda_z$,
then, of course, $\lambda_z=\phi_z^-(z)-\eta^-(z)$, and for every
central projection $z$ we chose $\phi_z^+$, $\phi_z^-$ such that
$\phi_z^+(z)=\eta^+(z)$, $\phi_z^-(z)=\eta^-(z)$. We can do it by
the following by
$$
\phi_z^+(x):=\phi_z^+(x)-\lambda_z,
\phi_z^-(x):=\phi_z^-(x)-\lambda_z, x\in z\mathcal{C}.
$$
We may take the following mapping
$$
\phi_+(za)=\phi_z^+(za), z\in Z(SC_{\#}(X,B(H)))\cap
P(SC_{\#}(X,B(H))), a\in \mathcal{C}.
$$

For every element $\sum_{i=1}^{m} x_ia_i$ from $\Re$  we suppose
$$
\phi_+(\sum_{i=1}^{m} x_ia_i)= \sum_{i=1}^{m} \phi_+(x_ia_i).
$$

Then $\phi_+$ is linear on $\Re$. Indeed, for arbitrary elements
$\sum_{i=1}^{m_1} x_ia_i$ and $\sum_{j=1}^{m_2} y_jb_j$ of $\Re$
we have
$$
\phi_+ (\sum_{i=1}^{m_1} a_iz_i+\sum_{j=1}^{m_2} y_jb_j)=\phi_+
[\sum_{i=1}^{m_1} \sum_{j=1}^{m_2} x_iy_j(a_i+b_j)]=
$$

$$
      \sum_{i=1}^{m_1} \sum_{j=1}^{m_2} \phi_+ (x_iy_ja_i+
      x_iy_jb_j)=\sum_{i=1}^{m_1}\sum_{j=1}^{m_2}\phi_+ (x_iy_ja_i)+
$$
$$
      \sum_{i=1}^{m_1} \sum_{j=1}^{m_2}\phi_+ (x_iy_jb_j)=
      \phi_+(\sum_{i=1}^{m_1} x_ia_i)+\phi_+(\sum_{i=1}^{m_2} y_jb_j).
$$
Thus $\phi_+ $ is linear on $\Re$. Similarly we may define the
linear mapping $\phi_-$ on $\Re$.

Note that from $\sum_{i=1}^{m_1} a_iz_i\geq 0$ it follows that
$a_iz_i\geq 0$ for every $i$. Hence, since $\phi_+(a_iz_i)\geq 0$
for each $i$ and $\phi_+$ is linear on $\Re$ we have
$\phi_+(\sum_{i=1}^{m_1} a_iz_i)\geq 0$. So, since
$\sum_{i=1}^{m_1} a_iz_i\geq 0$ is chosen arbitrarily in $\Re$ we
have $\phi_+(x)\geq 0$ for all $x\in\Re\cap C(X,B(H))_+$, where
$C(X,B(H))_+$ is the set of all positive elements of $C(X,B(H))$.

Let $\varepsilon >0$. Suppose that
$$
\Vert \sum_{i=1}^{m_1} x_ia_i-\sum_{j=1}^{m_2} y_jb_j \Vert
<\varepsilon.
$$

Since $\{x_iy_j\}$ is an orthogonal family of central projections
we have
$$
\Vert x_iy_j(a_i-b_j)\Vert<\varepsilon \eqno{(1.1)}
$$
for all indices $i$ and $j$.

Let us $x_iy_j(a_i-b_j)$ denote by $c_{ij}$ for every pair of
indices $i$, $j$. There exist pairwise orthogonal projections
$p^{ij}_1,p^{ij}_2,\dots,p^{ij}_k$ in $x_iy_j\mathcal{C}$ and
numbers $\beta^{ij}_1,\beta^{ij}_2,\dots, \beta^{ij}_k$ such that
$$
\Vert c_{ij}-\sum_{l=1}^k \beta^{ij}_lp^{ij}_l\Vert<\varepsilon .
\eqno{(1.2)}
$$
Since $\{x_iy_j\}$ is an orthogonal family of central projections
with sum $1$ we have
$$
\sum_{i=1}^{m_1} \sum_{j=1}^{m_2} x_iy_jc_{ij}<\varepsilon 1+
\sum_{i=1}^{m_1} \sum_{j=1}^{m_2}
[x_iy_j(\sum_{l=1}^k\beta^{ij}_lp^{ij}_l)],
$$
$$
-\sum_{i=1}^{m_1} \sum_{j=1}^{m_2} x_iy_jc_{ij}<\varepsilon 1-
\sum_{i=1}^{m_1} \sum_{j=1}^{m_2}
[x_iy_j(\sum_{l=1}^k\beta^{ij}_lp^{ij}_l)] \eqno{( 1.3}).
$$
The elements
$$
x_iy_j, x_iy_jc_{ij}, x_iy_j (\sum_{l=1}^k\beta^{ij}_lp^{ij}_l)
\,\,\text{for}\,\, \text{all}\,\,\text{indices}\,\,i,j,
$$
and
$$
\sum_{v=1}^{m_1} \sum_{w=1}^{m_2} x_vy_wc_{vw}, \sum_{v=1}^{m_1}
\sum_{w=1}^{m_2} [x_vy_w(\sum_{l=1}^k\beta^{vw}_lp^{vw}_l)]
$$
belong to $\Re$. Therefore the map $\phi_+$ is defined on these
elements. Let
$$
a^+=\phi_+(\sum_{i=1}^{m_1} \sum_{j=1}^{m_2} x_iy_jc_{ij}),
a^-=\phi_-(\sum_{i=1}^{m_1} \sum_{j=1}^{m_2} x_iy_jc_{ij}),
a=a^+-a^-,
$$
$$
b^+=\phi_+(\sum_{i=1}^{m_1} \sum_{j=1}^{m_2}[x_iy_j
(\sum_{l=1}^k\beta^{ij}_lp^{ij}_l)]), b^-=\phi_-(\sum_{i=1}^{m_1}
\sum_{j=1}^{m_2}[x_iy_j (\sum_{l=1}^k\beta^{ij}_lp^{ij}_l)]),
b=b^+-b^-.
$$
Then, since $\phi_+$ is linear and positive on $\Re$ we have
$$
a^+<\varepsilon \sum\phi_+(x_iy_j)+\sum_{i=1}^{m_1}
\sum_{j=1}^{m_2}\phi_+(x_iy_j (\sum_{l=1}^k\beta^{ij}_lp^{ij}_l))
$$
$$
=\varepsilon \phi_+(1)+b^+ \eqno{(1.4)}
$$
by  (1.3). Similarly
$$
-a^+<\varepsilon \phi_+(1)-b^+. \eqno{(1.5)}
$$

Similarly for $\phi_-$ we have
$$
a^-<\varepsilon \phi_-(1)+b^- \eqno{(1.4.1)}
$$
and
$$
-a^-<\varepsilon \phi_-(1)-b^-. \eqno{(1.5.1)}
$$
Hence
$$
a-b<\varepsilon (\phi_+(1)+\phi_-(1)), -\varepsilon
(\phi_+(1)+\phi_-(1))<a-b
$$
and
$$
\vert a-b\vert<\varepsilon (\phi_+(1)+\phi_-(1))<\varepsilon C_1,
$$
for some real number $C_1$ which is not depending on
$\varepsilon$.

Clearly the family $\{x_iy_j, x_iy_jp_l^{ij}\}_{ijl}$ of
projections is contained in some maximal commutative
$*$-subalgebra $A_\circ$ of $SC_{\#}(X,B(H))$.

The extension of $\phi$ to $A_\circ$ coincides with
$\delta\vert_{A_\circ}$ first by theorem 1.1 for every $i$ and $j$
the linear functional $\phi\vert_{x_iy_j\mathcal{C}}$ coincides on
$x_iy_j\mathcal{C}$ with the unique quasilinear functional
$\delta$, and second $\phi\vert_{P(A_\circ)}=\mu
\vert_{P(A_\circ)}$.

By (1.1) and (1.2)
$$
\Vert x_iy_j(\sum_{l=1}^k\beta^{ij}_lp^{ij}_l) \Vert\leq
2\varepsilon.
$$
At the same time  as it mentioned above $\phi$ is a continuous
linear functional on $A_\circ$. Therefore
$$
\vert\phi(\sum_{i=1}^{m_1} \sum_{j=1}^{m_2}
[x_iy_j(\sum_{l=1}^k\beta^{ij}_lp^{ij}_l)])\vert<2C_2\varepsilon,
$$
for some real number $C_2$ which is not depending on
$\varepsilon$. Hence by (1.4), (1.5), (1.4.1) and (1.5.1) we have
$$
\vert\phi(\sum_{i=1}^{m_1} \sum_{j=1}^{m_2}
x_iy_jc_{ij})\vert=\vert a\vert=\vert a-b+b\vert<\vert
a-b\vert+\vert b\vert <\varepsilon C_1+2C_2\varepsilon .
$$
and for every fundamental consequence $(a_n)\subset \Re$ the
numerical consequence $(\phi(a_n))$ is also fundamental.

Let $a$ be an arbitrary continuous operator-valued function from
$SC_{\#}(X,B(H))_{sa}$. For every natural number  $m$ and for each
point $x\in X$ there exists an open neighborhood $O_x\subseteq X$
of the point $x$ such that
$$
\Vert a(x)-a(y)\Vert <1/m
$$
for every $y\in O_x$. Fix $m$. Let $Q_x$ be the closure of the
neighborhood  $O_x$ for every point $x\in X$. Then for every point
$x\in X$ the set $Q_x$ is open. It is clear that  $\bigcup_{x\in
X} Q_x=X$. Thus the family $\{Q_x\}_{x\in X}$ is an open covering
of the compact $X$ and it can be chosen a finite covering, say
$\{Q_{x_k}\}_{k=1}^l$ in $\{Q_x\}_{x\in X}$. For every $k$ in
$\{1,2,\dots,l\}$ and $y\in Q_{x_k}$ we have
$$
\Vert a(x_k)-a(y)\Vert \leq 1/m.
$$
Without loss of the generality we admit that sets in
$\{Q_{x_k}\}_{k=1}^l$ are not pairwise crossed. For each $k$ in
$\{1,2,\dots,l\}$  the characteristic function
$\mathcal{X}_{Q_{x_k}}$ of the set $Q_{x_k}$  belongs to
$P(Z(A))$. Hence the element $\sum_{k=1}^l\mathcal{X}_{Q_{x_k}}
a(x_k)$  lies in $\Re$ and
$$
\Vert a-\sum_{k=1}^l\mathcal{X}_{Q_{x_k}} a(x_k)\Vert \leq 1/m.
$$
We continue the extension of $\phi$ on the uniform closure of
$\Re$ in $C(X,B(H)_{sa})$, which coincides with $C(X,B(H)_{sa})$
by the above conclusions. Take a consequence $(a_m)\subset\Re$
converging to the element $a$. Then, since $\phi$ is continuous on
$\Re$ there exists $\lim \phi(a_m)$. We admit $\phi(a)=\lim
\phi(a_m)$. It is clear that $\phi(a)$ does not depends on choice
of the consequence $(a_m)$ in $\Re$.

Let $p$ be a projection in $SC_{\#}(X,B(H))$, which is a
continuous operator-valued function. We prove that for $p$ the
element $\phi(a)$ is defined correctly. Let $\varepsilon$ be a
positive number. Then there exists an orthogonal family
$\{z_i\}_{i=1}^l$ of central projections such that $(\forall
i)z_i=\mathcal{X}_{Q_{x_i}}$ for some point $x_i\in X$ and
open-closed set $Q_{x_i}$ containing this point. In this case
$$
(\forall i)(\forall y\in Q_{x_i})\Vert p(x_i)-p(y)\Vert
<\varepsilon /2.
$$
Hence
$$
(\forall i)(\forall x,y\in Q_{x_i}) \Vert p(x)-p(y)\Vert
<\varepsilon.
$$
Fix  $i$ and $y$.  Then
$$
(\forall x\in Q_{x_i})-1\varepsilon<p(x)-p(y) <1 \varepsilon,
$$
i.e.
$$
-1\varepsilon+p(y)<p(x) <1 \varepsilon+p(y).
$$
Hence
$$
-1\varepsilon+p(y)<\inf_{x\in Q_{x_i}} p(x) <1 \varepsilon+p(y),
$$
i.e.
$$
-1\varepsilon<\inf_{x\in Q_{x_i}} p(x)-p(y) <1 \varepsilon.
$$
Therefore
$$
(\forall y\in Q_{x_i}) \Vert \inf_{x\in Q_{x_i}} p(x)-p(y) \Vert<
\varepsilon.
$$
Let $(\forall i)a_i=\inf_{x\in Q_{x_i}} p(x)$. Then
$$
\Vert p-\sum z_ia_i\Vert<\varepsilon, \sum z_ia_i\in \Re
$$
and
$$
p\geq \sum z_ia_i.
$$
Besides, since $(\forall i) (\forall x\in Q_{x_i}) p(x)\geq 0$ we
have $\sum z_ia_i\geq 0$ ($\sum z_ia_i$ is a projection). Hence
$$
0=(1-p)p(1-p)\geq (1-p)(\sum z_ia_i)(1-p)\geq 0
$$
and
$$
p(\sum z_ia_i)=\sum z_ia_i.
$$
Hence $p$, $\sum z_ia_i$ mutually commute in $\mathcal{A}$.
Clearly the set
$$
\{p, \sum z_ia_i, z_1a_1, z_2a_2, \dots, z_la_l\}
$$
is a set of pairwise commutative elements in $C(X,B(H))_{sa}$. For
the maximal commutative $*$-subalgebra $A_\circ$ of
$SC_{\#}(X,B(H))$, containing this set the value at $\sum z_ia_i$
of $\delta$ on $A_\circ$, to which the measure $\mu$ is extended,
coincides with the value of the function $\phi$ at this element
$\sum z_ia_i$. Indeed, by theorem 1.1 for every $i$ the linear
functional $\phi\vert_{z_i\mathcal{C}}$ coincides on
$z_i\mathcal{C}$ with $\delta$. Therefore $(\forall
i)\delta(z_ia_i)=\phi(z_ia_i)$. Hence $\delta(\sum
z_ia_i)=\phi(\sum z_ia_i)$. At the same time, since the measure
$\mu$ is uniquely extended to $\delta$ on $SC_{\#}(X,B(H))$
$\phi(\sum z_ia_i)$ does not depend on choice of $A_\circ$. Since
$\delta$ is continuous on $A_\circ$ and $\Vert p-\sum z_ia_i\Vert<
\varepsilon$ it follows that
$$
\vert \mu(p)-\phi(\sum z_ia_i)\vert<\varepsilon.
$$

Thus, there exists a consequence $(a_m)$ in $\Re$ uniformly
converging to $p$ such that $p$ and $a_m$ are mutually commutes
for each $m$. In this case $\lim \phi(a_m)$ exists and
$\mu(p)=\lim \phi(a_m)$. So the extension of $\phi\vert_{\Re}$ on
$C(X,B(H))_{sa}$ coincides on projections with the unique
quasi-linear functional $\delta$ defined on $SC_{\#}(X,B(H))$ to
which the measure $\mu$ is extended. Then $\phi$ is linear on
$C(X,B(H))_{sa}$. Let $a,b\in C(X,B(H))_{sa}$ such that $\Vert a-b
\Vert <\varepsilon$. Then there exist $c$ and $d$ in $\Re$ such
that $\Vert a-c \Vert <\varepsilon $ and $\Vert b-d\Vert
<\varepsilon $. Hence
$$
\vert \phi (a-b)\vert = \vert\phi (a-c+c-d+d-b)\vert\leq
$$
$$
\vert \phi (a-c)\vert +\vert \phi (c-d)\vert +\vert \phi
(d-b)\vert <C\varepsilon,
$$
where $C$ is a constant that does not depend on $\varepsilon$.

Thus, $\phi $ is a continuous linear functional on the subspace
$C(X,B(H)_{sa})$ of $SC_{\#}(X,B(H))_{sa}$. Hence the map
$$
\rho(a+ib)=\phi(a)+i\phi(b), a, b\in C(X,B(H)_{sa})
$$
is a linear functional on $C(X,B(H))$, which is a unique extension
of $\mu$. $\triangleright$

{\bf Corollary 1.3.} The set of all linear combinations of all
finite families of orthogonal projections in $C(X,B(H))$ is
uniformly dense in $C(X,B(H))$.

{\bf Remark.} Note that theorem 1.1 was proved in the case of
bounded signed measure by Matvejchuk M.S in \cite{MMS} for purely
infinite AW$^*$-algebras and finite AW$^*$-algebras with a
faithful normal centrevalued trace.

{\bf Lemma 1.4.} {\it Let $\mathcal{A}$ be a AW$^*$-algebra,
$\Delta :\mathcal{A}\to \mathcal{A}$ be a 2-local derivation, and
let $e$, $f$ be mutually orthogonal projections in $\mathcal{A}$.
Then

1) $\Delta(e+f)=\Delta(e)+\Delta(f)$,

2) if $\lambda$, $\mu$  are arbitrary complex numbers, then
$\Delta(\lambda e+\mu f)=\Delta(\lambda e)+\Delta(\mu f)$.}

{\bf Proof.} A proof of 1): by the definition there exist $a$,
$b\in \mathcal{A}$ such that
$$
\Delta(e+f)=a(e+f)-(e+f)a,  \Delta(e)=ae-ea,
$$
$$
\Delta(e+f)=b(e+f)-(e+f)b,  \Delta(f)=bf-fb.
$$
Then
$$
(e+f)a(e+f)^\perp=(e+f)b(e+f)^\perp,
$$
$$
(e+f)^\perp a(e+f)= (e+f)^\perp b(e+f).
$$
Hence
$$
fa(e+f)^\perp=fb(e+f)^\perp, (e+f)^\perp af=(e+f)^\perp bf.
$$

Now we must show $eaf=ebf$, $fae=fbe$. Indeed, there exists $d\in
M$ such that
$$
\Delta(e)=de-ed,  \Delta(f)=df-fd.
$$
Hence
$$
fae=fde, eaf=edf, ebf=edf, fbe=fde.
$$
Therefore
$$
eaf=ebf, fae=fbe
$$
and
$$
\Delta(e+f)=a(e+f)-(e+f)a=
$$
$$
ea(e+f)+fa(e+f)+(e+f)^\perp a(e+f)-
$$
$$
(e+f)ae-(e+f)af-(e+f)a(e+f)^\perp=
$$
$$
eae+fae+(e+f)^\perp ae-eae-eaf-ea(e+f)^\perp+
$$
$$
eaf+faf+(e+f)^\perp af-fae-faf-fa(e+f)^\perp=
$$
$$
eae+fae+(e+f)^\perp ae-eae-eaf-ea(e+f)^\perp+
$$
$$
ebf+fbf+(e+f)^\perp bf-fbe-fbf-fb(e+f)^\perp=
$$
$$
ae-ea+bf-fb=\Delta(e)+\Delta(f).
$$

A proof of 2): by the definition there exist $a$, $b\in
\mathcal{A}$ such that
$$
\Delta(\lambda e+\mu f)=a(\lambda e+\mu f)-(\lambda e+\mu f)a,
\Delta(e)=ae-ea,
$$
$$
\Delta(\lambda e+\mu f)=b(\lambda e+\mu f)-(\lambda e+\mu f)b,
\Delta(f)=bf-fb.
$$
Then
$$
(\lambda e+\mu f)a(e+f)^\perp=(\lambda e+\mu f)b(e+f)^\perp,
$$
$$
(e+f)^\perp a(\lambda e+\mu f)= (e+f)^\perp b(\lambda e+\mu f),
$$
$$
eaf=ebf, fae=fbe.
$$
Hence
$$
fa(e+f)^\perp =fb(e+f)^\perp, (e+f)^\perp af=(e+f)^\perp bf
$$
and
$$
\Delta (\lambda e+\mu f)=a(\lambda e+\mu f)-(\lambda e+\mu f)a=
$$
$$
ea(\lambda e+\mu f)+fa(\lambda e+\mu f)+(e+f)^\perp a(\lambda
e+\mu f)-
$$
$$
(\lambda e+\mu f)ae-(\lambda e+\mu f)af-(\lambda e+\mu
f)a(e+f)^\perp=
$$
$$
ea\lambda e+fa\lambda e+(e+f)^\perp a\lambda e-\lambda eae-\lambda
eaf- \lambda ea(e+f)^\perp +
$$
$$
ea\mu f+fa\mu f+(e+f)^\perp a\mu f-\mu fae-\mu faf-\mu
fa(e+f)^\perp =
$$
$$
ea\lambda e+fa\lambda e+(e+f)^\perp a\lambda e-\lambda eae-\lambda
eaf- \lambda ea(e+f)^\perp +
$$
$$
eb\mu f+fb\mu f+(e+f)^\perp b\mu f-\mu fbe-\mu fbf-\mu
fb(e+f)^\perp=
$$
$$
\lambda (ae-ea)+\mu (bf-fb)=\lambda  \Delta(e)+\mu \Delta(f)=
$$
$$
\Delta(\lambda e)+\Delta(\mu f).
$$
$\triangleright$

\medskip

{\bf Lemma 1.5.} {\it Let , $\Delta :\mathcal{A}\to \mathcal{A}$
be a 2-local derivation on a AW$^*$-algebra $\mathcal{A}$. Then
$$
\Delta(\lambda_1e_1+\lambda_2e_2+...+\lambda_me_m)=
$$
$$
\Delta(\lambda_1e_1)+\Delta(\lambda_2e_2)+...+\Delta(\lambda_me_m).
$$
for any family of pairwise orthogonal projections $e_1$, $e_2$,
$\dots$, $e_m$ in $\mathcal{A}$ and complex numbers $\lambda_1$,
$\lambda_2$, ...,$\lambda_m$.}

{\bf Proof.} It is clear that
$$
\Delta(e_1+e_2+...+e_m)=\Delta (e_1)+\Delta (e_2+...+e_m)=
$$
$$
\Delta (e_1)+\Delta (e_2)+\Delta (e_3+...+e_m)=...=
$$
$$
\Delta (e_1)+\Delta (e_2)+...+\Delta (e_m)
$$
by 1) of lemma 1.4.

Using the induction we prove that
$$
\Delta(\lambda_1e_1+\lambda_2e_2+...+\lambda_me_m)=
$$
$$
\Delta(\lambda_1e_1)+\Delta(\lambda_2e_2)+...+\Delta(\lambda_me_m).
$$
The case $m=1$ is obvious. The case $m=2$ follows by lemma 1.4.
Suppose that
$$
\Delta(\lambda_1e_1+\lambda_2e_2+...+\lambda_{m-1}e_{m-1})=
$$
$$
\Delta(\lambda_1e_1)+\Delta(\lambda_2e_2)+...+\Delta(\lambda_{m-1}e_{m-1}).
$$
Let $x=\lambda_1e_1+\lambda_2e_2+...+\lambda_{m-1}e_{m-1}$. Then
there exist $a$, $b\in \mathcal{A}$ such that
$$
\Delta(\lambda_me_m+x)=
$$
$$
a(\lambda_me_m+x)-(\lambda_me_m+x)a, \Delta(e_m)=ae_m-e_ma,
$$
$$
\Delta(\lambda_me_m+x)=b(\lambda_me_m+x)-(\lambda_me_m+x)b,
$$
$$
\Delta(x)=bx-xb.
$$
Then
$$
(\lambda_me_m+x)a(e_m+f)^\perp=(\lambda_me_m+x)b(e_m+f)^\perp,
$$
$$
(e_m+f)^\perp a(\lambda_me_m+x)=(e_m+f)^\perp b(\lambda_me_m+x),
$$
where $f=e_1+e_2+...+e_{m-1}$. Hence
$$
xa(e_m+f)^\perp=xb(e_m+f)^\perp,
$$
$$
(e_m+f)^\perp ax=(e_m+f)^\perp bx, \,\,\,\,\,\,\,\,\,\, (1.6)
$$
$$
e_ma(e_m+f)^\perp =e_mb(e_m+f)^\perp,
$$
$$
(e_m+f)^\perp ae_m=(e_m+f)^\perp be_m.
$$
Let us prove that
$$
e_max=e_mbx, xae_m=xbe_m. \,\,\,\,\,\,\,\,\,\, (1.7)
$$

Indeed, there exists $d\in \mathcal{A}$ such that
$$
\Delta (e_m)=de_m-e_md,  \Delta(x)=dx-xd.
$$
Hence
$$
xae_m=xde_m, e_max=e_mdx,
$$
$$
e_mbx=e_mdx, xbe_m=xde_m.
$$
Therefore
$$
e_max=e_mbx, xae_m=xbe_m.
$$
Hence
$$
e_maf=e_mbf, fae_m=fbe_m.
$$
Also we have
$$
(e_m+f+(e_m+f)^\perp)a(\lambda_me_m+x)-
$$
$$
(\lambda_me_m+x)a(e_m+f+(e_m+f)^\perp)=
$$
$$
(e_m+f+(e_m+f)^\perp)b(\lambda_me_m+x)-
$$
$$
(\lambda_me_m+x)b(e_m+f+(e_m+f)^\perp),
$$
and
$$
fax-xaf=fbx-xbf. \,\,\,\,\,\,\,\,\,\,  (1.8)
$$
Therefore
$$
\Delta(\lambda_me_m+x)=a(\lambda_me_m+x)-(\lambda_me_m+x)a=
$$
$$
e_ma(\lambda_me_m+x)+fa(\lambda_me_m+x)+(e_m+f)^\perp
a(\lambda_me_m+x)-
$$
$$
(\lambda_me_m+x)ae_m-(\lambda_me_m+x)af-(\lambda_me_m+x)a(e_m+f)^\perp=
$$
$$
e_ma\lambda_me_m+fa\lambda_me_m+(e_m+f) a\lambda_me_m-
$$
$$
\lambda_me_mae_m-\lambda_me_maf-\lambda_me_ma(e_m+f)^\perp+
$$
$$
e_max+fax+(e_m+f)^\perp ax-xae_m-xaf-xa(e_m+f)^\perp=
$$
$$
e_ma\lambda_me_m+fa\lambda_me_m+(e_m+f)^\perp a\lambda_me_m-
$$
$$
\lambda_me_mae_m-\lambda_me_maf-\lambda_me_ma(e_m+f)^\perp+
$$
$$
e_mbx+fbx+(e_m+f)^\perp bx-xbe_m-xbf-xb(e_m+f)^\perp=
$$
$$
\lambda_m(ae_m-e_ma)+(bx-xb)=\lambda_m \Delta(e_m)+ \Delta(x)=
$$
$$
\Delta(\lambda_me_m)+\Delta(x)
$$
by (1.6), (1.7) and (1.8).

Hence by induction we obtain that
$$
\Delta(\lambda_1e_1+\lambda_2e_2+...+\lambda_me_m)=
$$
$$
\Delta(\lambda_1e_1)+\Delta(\lambda_2e_2)+...+\Delta(\lambda_me_m).
$$
$\triangleright$

{\bf Lemma 1.6.} {\it Let $\mathcal{A}$ be a AW$^*$-algebra,
$\Delta :\mathcal{A}\to \mathcal{A}$ be a 2-local derivation, and
let $\mathcal{A}_o$ be a maximal abelian $*$-subalgebra of
$\mathcal{A}$. Consider two linear combinations,
$\lambda_1e_1+\lambda_2e_2+...+\lambda_me_m$,
$\mu_1f_1+\mu_2f_2+...+\mu_kf_k$ defined by the sets
$\{e_1,e_2,...,e_m\}$, $\{f_1,f_2,...,f_k\}$ of orthogonal
projections in $\mathcal{A}_o$ respectively. Then
$$
\Delta((\lambda_1e_1+\lambda_2e_2+...+\lambda_me_m)+(\mu_1f_1+
\mu_2f_2+...+\mu_kf_k))=
$$
$$
\Delta(\lambda_1e_1+\lambda_2e_2+...+\lambda_me_m)+\Delta(\mu_1f_1+
\mu_2f_2+...+\mu_kf_k).
$$
}

{\bf Proof.} It clear that the union of the families
$$
\{e_i-e_if\}_{i=1,...,m}, \{e_if_j\}_{i=1,...,m,j=1,...,k},
\{f_j-f_je\}_{j=1,...,k}
$$
is a set of orthogonal projections in $M_o$, where
$e=e_1+e_2+...+e_m$, $f=f_1+f_2+...+f_k$. Hence
$$
\Delta(\lambda_1e_1+\lambda_2e_2+...+\lambda_me_m+\mu_1f_1+\mu_2f_2+...+\mu_kf_k)=
$$
$$
\Delta(\sum_{i=1}^m
\lambda_i(e_i-e_if)+\sum_{i=1}^m\sum_{j=1}^k(\lambda_i+\mu_j)e_if_j+
\sum_{i=1}^k \mu_i(f_i-f_ie))=
$$
$$
\sum_{i=1}^m\lambda_i\Delta(e_i-e_if)+\sum_{i=1}^m\sum_{j=1}^k(\lambda_i+\mu_j)
\Delta(e_if_j)+\sum_{i=1}^k\mu_i \Delta(f_i-f_ie))=
$$
$$
[\sum_{i=1}^m \lambda_i \Delta(e_i-e_if)+ \sum_{i=1}^m\sum_{j=1}^k
i \Delta(e_if_j)]+
$$
$$
[\sum_{i=1}^k \mu_i \Delta(f_i-f_ie))+\sum_{i=1}^m\sum_{j=1}^k
\mu_j \Delta(e_if_j)]=
$$
$$
\Delta(\lambda_1e_1+\lambda_2e_2+...+\lambda_me_m)+\Delta(\mu_1f_1+
\mu_2f_2+...+\mu_kf_k)
$$
by lemma 1.5. $\triangleright$

In the following lemmas $\mathcal{P}(C(X,B(H)))$ denotes the
lattice of all projections in $C(X,B(H))$ and
$Ln_o(\mathcal{P}(C(X,B(H))))$ denotes the set of all finite
linear combinations of orthogonal projections in $C(X,B(H))$.

{\bf Lemma 1.7.} {\it For every pair of elements $x$, $y\in
Ln_o(\mathcal{P}(C(X,B(H))))$
$$
\Delta(x+y)=\Delta(x)+\Delta(y).
$$
}

{\bf Proof.} Firstly, let us show that for each $f\in
SC_{\#}(X,B(H))^*$ the restriction
$f\circ\Delta\vert_{\mathcal{P}(SC_{\#}(X,B(H)))}$ of the
superposition $f\circ\Delta(x)=f(\Delta(x))$, $x\in
SC_{\#}(X,B(H))$, onto the lattice $\mathcal{P}(SC_{\#}(X,B(H)))$
is a bounded signed measure, where $SC_{\#}(X,B(H))^*$ is the
space of linear functionals on $SC_{\#}(X,B(H))$. Let $e_1$, $e_2$
be orthogonal projections in $SC_{\#}(X,B(H))$. By 1) of lemma 1.4
we obtain that
$$
f\circ\Delta(e_1+e_2)=f(\Delta(e_1)+\Delta(e_2))=f(\Delta(e_1))+f(\Delta(e_2))=
$$
$$
f\circ\Delta(e_1)+f\circ\Delta(e_2),
$$
i.e. $f\circ\Delta$ is a signed measure.

By theorem 1.2 (Gleason Theorem) for signed measures there exists
a unique bounded linear functional $\tilde{f}$ on $C(X,B(H))$ such
that
$$
\tilde{f}\vert_{\mathcal{P}(C(X,B(H)))}=f\circ\Delta\vert_{\mathcal{P}(C(X,B(H)))}.
$$
Let us show that
$\tilde{f}\vert_{Ln_o(\mathcal{P}(C(X,B(H))))}=f\circ\Delta\vert_{Ln_o(\mathcal{P}(C(X,B(H))))}$.
Indeed, let $\mathcal{A}_o$ be a maximal commutative
$*$-subalgebra of $C(X,B(H))$. Then by lemma 1.6 $\Delta$ is
linear on $\mathcal{A}_o\cap Ln_o(\mathcal{P}(C(X,B(H))))$, and
therefore $f\circ\Delta\vert_{\mathcal{A}_o\cap
Ln_o(\mathcal{P}(C(X,B(H))))}$ is a bounded linear functional
which is an extension of the measure
$f\circ\Delta\vert_{\mathcal{P}(C(X,B(H)))}$. By the uniqueness of
the extension we have $\tilde{f}(x)=f\circ\Delta(x)$, $x\in
\mathcal{A}_o\cap Ln_o(\mathcal{P}(C(X,B(H))))$.

So for all $f\in SC_{\#}(X,B(H))^*$ we have
$$
f(\Delta(x+y))=f(\Delta(x))+f(\Delta(y))=f(\Delta(x)+\Delta(y)),
$$
i.e. $f(\Delta(x+y)-\Delta(x)-\Delta(y))=0$ for all $f\in
SC_{\#}(X,B(H))^*$. Since $SC_{\#}(X,B(H))^*$ separates points of
$\mathcal{A}$ it follows that $\Delta(x+y)-\Delta(x)-\Delta(y)=0$
for all $x$, $y\in Ln_o(\mathcal{P}(C(X,B(H))))$. The proof is
complete. $\triangleright$

{\bf Lemma 1.8.} {\it There exists an element $a\in
C(X,B(H))^{**}$ such that $\Delta(x)=D_a(x)=ax-xa$ for all $x\in
Ln_o(\mathcal{P}(C(X,B(H))))$.}

{\bf Proof.} We have the set $Ln_o(\mathcal{P}(C(X,B(H))))$ is
uniformly dense in $C(X,B(H))$ by corollary 1.3 and
$C(X,B(H)_{sa})$ is weakly dense in $C(X,B(H))^{**}_{sa}$.
Therefore $Ln_o(\mathcal{P}(C(X,B(H))))$ is weakly dense in
$C(X,B(H))^{**}_{sa}$. Let $x$ be an arbitrary element in
$C(X,B(H))^{**}_{sa}$ and $(x_n)$ be a sequence in
$Ln_o(\mathcal{P}(C(X,B(H))))$ weakly converging to $x$. Then by
the proof of lemma 1.7 for any $f\in C(X,B(H))^{***}$ the sequence
$(f\circ \Delta (x_n))$ is a fundamental sequence of complex
numbers. Hence the sequence $\Delta (x_n)$ is also a fundamental
sequence in the weak topology and weakly converges to some element
$y\in C(X,B(H))^{**}$. Let $\tilde{\Delta}(x)=y$. Consider the
weak extension $\tilde{\Delta}$ of
$\Delta\vert_{Ln_o(\mathcal{P}(C(X,B(H))))}$ on
$C(X,B(H))^{**}_{sa}$. Then this extension $\tilde{\Delta}$ is
additive on $C(X,B(H))^{**}_{sa}$ by lemma 1.7. Taking into
account the homogeneity of $\Delta$ we obtain that
$$
\tilde{\Delta}(x^2)=\tilde{\Delta}(x)x+x\tilde{\Delta}(x), x\in
C(X,B(H))^{**}_{sa}
$$
since multiplication is separately weakly continuous in
$C(X,B(H))^{**}$. Consider the extension $\hat{\Delta}$ of
$\tilde{\Delta}\vert_{C(X,B(H))^{**}_{sa}}$ on $C(X,B(H))^{**}$
defined by:
$$
\hat{\Delta}(x_1+ix_2) =\tilde{\Delta}(x_1)+i\tilde{\Delta}(x_2),
x_1, x_2\in C(X,B(H))^{**}_{sa}.
$$
By the definition $\tilde{\Delta}$ is a Jordan derivation on
$C(X,B(H))^{**}$. As we mentioned above by \cite[Theorem 1]{Bre}
any Jordan derivation on a semiprime algebra is a derivation.
Since $C(X,B(H))^{**}$ is semi-prime $\tilde{\Delta}$ is a
derivation on $C(X,B(H))^{**}$. It is known \cite{S} that any
derivation $D$ on $C(X,B(H))^{**}$ is an inner derivation, that is
there exists an element $a\in C(X,B(H))^{**}$ such that
$D(x)=ax-xa$ for all $x\in C(X,B(H))^{**}$. Therefore there exists
an element $a\in C(X,B(H))^{**}$ such that
$$
\tilde{\Delta}(x)=ax-xa
$$
for all $x\in C(X,B(H))^{**}$. In particular,
$\Delta(x)=D_a(x)=ax-xa$ for all $x\in
Ln_o(\mathcal{P}(C(X,B(H))))$. The proof is complete.
$\triangleright$

\medskip

\section{2-local derivations on AW$^*$-algebras of type I$_n$}

\medskip

We take the AW$^*$-algebra of type I$_n$ $SC_{\#}(X,B(H))$, the
C$^*$-subalgebra $C(X,B(H))$ of all continuous operator-valued
functions on $X$ and the system of matrix units $\{e_{ij}\}$ from
section 1.

\medskip

{\bf Lemma 2.1.} {\it Let $a$ be the element from lemma 1.8. Then
for every $i$ and $j$
$$
e_ia(ij)e_j=e_iae_j, e_ja(ij)e_i=e_jae_i,
$$
$$
e_ia(ij)e_{ij}-e_{ij}a(ij)e_j=e_iae_{ij}-e_{ij}ae_j,
$$
where $a(ij)\in SC_{\#}(X,B(H))$ is an element such that
$$
\bigtriangleup(e_{ij})=a(ij)e_{ij}-e_{ij}a(ij).
$$
}

{\it Proof.} There exists an element $d$ in $SC_{\#}(X,B(H))$ such
that
$$
\bigtriangleup(e_i)=de_i-e_id,
\bigtriangleup(e_{ij})=de_{ij}-e_{ij}d.
$$
Hence
$$
de_{ij}-e_{ij}d=a(ij)e_{ij}-e_{ij}a(ij)
$$
and
$$
e_ide_j=e_ia(ij)e_j, e_jde_i=e_ja(ij)e_i.
$$
At the same time
$$
de_i-e_id=ae_i-e_ia
$$
by lemma 1.8 and
$$
e_ide_j=e_iae_j, e_jde_i=e_jae_i.
$$
Therefore
$$
e_iae_j=e_ia(ij)e_j, e_jae_i=e_ja(ij)e_i.
$$

Now we have $e_{ij}+e_{ji}\in Ln_o(\mathcal{P}(C(X,B(H))))$.
Therefore
$$
\bigtriangleup(e_{ij}+e_{ji})=a(e_{ij}+e_{ji})-(e_{ij}+e_{ji})a
$$
by lemma 1.8. By the definition of 2-local derivation there exists
$d\in SC_{\#}(X,B(H))$ such that
$$
\bigtriangleup(e_{ij})=de_{ij}-e_{ij}d,
$$
$$
\bigtriangleup(e_{ij}+e_{ji})=d(e_{ij}+e_{ji})-(e_{ij}+e_{ji})d.
$$
Hence
$$
d(e_{ij}+e_{ji})-(e_{ij}+e_{ji})d=a(e_{ij}+e_{ji})-(e_{ij}+e_{ji})a
$$
and
$$
e_ide_{ij}-e_{ij}de_j=e_iae_{ij}-e_{ij}ae_j.
$$
At the same time, since
$$
de_{ij}-e_{ij}d=a(ij)e_{ij}-e_{ij}a(ij)
$$
we have
$$
e_ide_{ij}-e_{ij}de_j=e_ia(ij)e_{ij}-e_{ij}a(ij)e_j.
$$
Therefore
$$
e_ia(ij)e_{ij}-e_{ij}a(ij)e_j=e_iae_{ij}-e_{ij}ae_j.
$$
$\triangleright$

{\bf Lemma 2.2.} {\it Let $a$ be an element from lemma 1.8. Then
for any pair $i$, $j$ of different indices the following equality
holds
$$
\bigtriangleup(e_{ij})=ae_{ij}-e_{ij}a.    \,\,\,\,\,\,\,(2.1)
$$
}

{\it Proof.} Let $k$ be an arbitrary index different from $i$, $j$
and let $a(ij,ik)\in SC_{\#}(X,B(H))$ be an element such that
$$
\bigtriangleup(e_{ik})=a(ij,ik)e_{ik}-e_{ik}a(ij,ik) \,\,
\text{and}\,\,
\bigtriangleup(e_{ij})=a(ij,ik)e_{ij}-e_{ij}a(ij,ik).
$$
Then
$$
e_{kk}\bigtriangleup(e_{ij})e_{jj}=e_{kk}(a(ij,ik)e_{ij}-e_{ij}a(ij,ik))e_{jj}=
$$
$$
e_{kk}a(ij,ik)e_{ij}-0=e_{kk}a(ik)e_{ij}-e_{kk}e_{ij}ae_{jj}=
$$
$$
e_{kk}a_{ki}e_{ij}-e_{kk}e_{ij}ae_{jj}=
e_{kk}ae_{ij}-e_{kk}e_{ij}ae_{jj}=
$$
$$
e_{kk}(ae_{ij}-e_{ij}a)e_{jj}
$$
by lemma 2.1.

Similarly,
$$
e_{kk}\bigtriangleup(e_{ij})e_{ii}=e_{kk}(a(ij,ik)e_{ij}-e_{ij}a(ij,ik))e_{ii}=
$$
$$
e_{kk}a(ij,ik)e_{ij}e_{ii}-0=0-0=
e_{kk}ae_{ij}e_{ii}-e_{kk}e_{ij}ae_{ii}=
$$
$$
e_{kk}(ae_{ij}-e_{ij}a)e_{ii}.
$$

Let $a(ij,kj)\in SC_{\#}(X,B(H))$ be an element such that
$$
\bigtriangleup(e_{kj})=a(ij,kj)e_{kj}-e_{kj}a(ij,kj)  \,\,
\text{and}\,\,
\bigtriangleup(e_{ij})=a(ij,kj)e_{ij}-e_{ij}a(ij,kj).
$$

Then
$$
e_{ii}\bigtriangleup(e_{ij})e_{kk}=e_{ii}(a(ij,kj)e_{ij}-e_{ij}a(ij,kj))e_{kk}=
$$
$$
0-e_{ij}a(ij,kj)e_{kk}=0-e_{ij}a(kj)e_{kk}=0-e_{ij}a_{jk}e_{kk}=
$$
$$
e_{ii}ae_{ij}e_{kk}-e_{ij}ae_{kk}=
$$
$$
e_{ii}(ae_{ij}-e_{ij}a)e_{kk}
$$
by lemma 2.1.

Also we have
$$
e_{jj}\bigtriangleup(e_{ij})e_{kk}=e_{jj}(a(ij,kj)e_{ij}-e_{ij}a(ij,kj))e_{kk}=
$$
$$
0-0=e_{jj}\{a(ij)\}_{i\neq
j}e_{ij}e_{kk}-e_{jj}e_{ij}\{a(ij)\}_{i\neq j}e_{kk}=
$$
$$
e_{jj}(ae_{ij}-e_{ij}a)e_{kk},
$$

$$
e_{ii}\bigtriangleup(e_{ij})e_{ii}=e_{ii}(a(ij)e_{ij}-e_{ij}a(ij))e_{ii}=
$$
$$
0-e_{ij}a(ij)e_{ii}=0-e_{ij}a(ij)e_{ii}=0-e_{ij}a_{ji}e_{ii}=
$$
$$
e_{ii}ae_{ij}e_{ii}-e_{ij}ae_{ii}=
$$
$$
e_{ii}(ae_{ij}-e_{ij}a)e_{ii}
$$
by lemma 2.1.

$$
e_{jj}\bigtriangleup(e_{ij})e_{jj}=e_{jj}(a(ij)e_{ij}-e_{ij}a(ij))e_{jj}=
$$
$$
e_{jj}a(ij)e_{ij}-0=e_{jj}a_{ji}e_{ij}-0=
$$
$$
e_{jj}ae_{ij}-e_{jj}e_{ij}ae_{jj}=
$$
$$
e_{jj}(ae_{ij}-e_{ij}a)e_{jj}
$$
by lemma 2.1.

$$
e_{ii}\bigtriangleup(e_{ij})e_{jj}=e_{ii}(a(ij)e_{ij}-e_{ij}a(ij))e_{jj}=
$$
$$
e_{ii}a(ij)e_{ij}-e_{ij}a(ij)e_{jj}=e_{ii}ae_{ij}-e_{ij}ae_{jj},
$$

$$
e_{jj}\bigtriangleup(e_{ij})e_{ii}=e_{jj}(a(ij)e_{ij}-e_{ij}a(ij))e_{ii}=0=
$$
$$
e_{jj}(ae_{ij}-e_{ij}a)e_{ii},
$$
by lemma 2.1.

Therefore for all indices $\alpha$ and $\beta$ we have
$$
e_{\alpha\alpha}\bigtriangleup(e_{ij})e_{\beta\beta}=e_{\alpha\alpha}(ae_{ij}-e_{ij}a)e_{\beta\beta}.
$$
Hence the equality (2.1) holds. $\triangleright$

\medskip

{\bf Theorem 2.3.} {\it There exists an element $d\in
SC_{\#}(X,B(H))$ such that $\bigtriangleup(x)=dx-xd$ for every
$x\in SC_{\#}(X,B(H))$ and $\bigtriangleup$ is a derivation on
$SC_{\#}(X,B(H))$.}

{\it Proof.} Let $a$ be an element from lemma 1.8 and $d(ij)\in
SC_{\#}(X,B(H))$ be an element such that
$$
\bigtriangleup(e_{ij})=d(ij)e_{ij}-e_{ij}d(ij) \,\, \text{and}\,\,
\bigtriangleup(x)=d(ij)x-xd(ij)
$$
and $i\neq j$. Then
$$
d(ij)e_{ij}-e_{ij}d(ij)=ae_{ij}-e_{ij}a
$$
for all $i$, $j$ by lemma 2.2 and
$$
(1-e_{ii})d(ij)e_{ii}=(1-e_{ii})ae_{ii},
e_{jj}d(ij)(1-e_{jj})=e_{jj}a(1-e_{jj}), \,\,\,\,\,\,\,\,\,(2.2)
$$
$$
e_{ii}d(ij)e_{ij}-e_{ij}d(ij)e_{jj}=e_{ii}ae_{ij}-e_{ij}ae_{jj}.\,\,\,\,\,\,\,\,\,(2.3)
$$
for all different $i$ and $j$.

Hence by (2.2), (2.3) we have
$$
e_{jj}\bigtriangleup(x)e_{ii}=e_{jj}(d(ij)x-xd(ij))e_{ii}=
$$
$$
e_{jj}d(ij)(1-e_{jj})xe_{ii}+
e_{jj}d(ij)e_{jj}xe_{ii}-e_{jj}x(1-e_{ii})d(ij)e_{ii}-e_{jj}xe_{ii}d(ij)e_{ii}=
$$
$$
e_{jj}a(1-e_{jj})xe_{ii}-e_{jj}x(1-e_{ii})ae_{ii}+
e_{jj}d(ij)e_{jj}xe_{ii}-e_{jj}xe_{ii}d(ij)e_{ii}=
$$
$$
e_{jj}a(1-e_{jj})xe_{ii}-e_{jj}x(1-e_{ii})ae_{ii}+
e_{jj}ae_{jj}xe_{ii}-e_{jj}xe_{ii}ae_{ii}=
$$
$$
e_{jj}(ax-xa)e_{ii}
$$
for all different $i$ and $j$.

Let $d(ii)\in SC_{\#}(X,B(H))$ be an element such that
$$
\bigtriangleup(e_{ii})=d(ii)e_{ii}-e_{ii}d(ii) \,\, \text{and}\,\,
\bigtriangleup(x)=d(ii)x-xd(ii)
$$
for each $i$. Then
$$
d(ii)e_{ii}-e_{ii}d(ii)=ae_{ii}-e_{ii}a
$$
by lemma 1.8 for all $i$  and
$$
(1-e_{ii})d(ii)e_{ii}=(1-e_{ii})ae_{ii},
e_{ii}d(ii)(1-e_{ii})=e_{ii}a(1-e_{ii}), \,\,\,\,\,\,\,\,\,(2.4)
$$
$$
e_{ii}d(ii)e_{ii}-e_{ii}d(ii)e_{ii}=e_{ii}ae_{ii}-e_{ii}ae_{ii}=0.\,\,\,\,\,\,\,\,\,(2.5)
$$
for every $i$.

Also by (2.4), (2.5) we have
$$
e_{ii}\bigtriangleup(x)e_{ii}=e_{ii}(d(ii)x-xd(ii))e_{ii}=
$$
$$
e_{ii}d(ii)(1-e_{ii})xe_{ii}+
e_{ii}d(ii)e_{ii}xe_{ii}-e_{ii}x(1-e_{ii})d(ii)e_{ii}-e_{ii}xe_{ii}d(ii)e_{ii}=
$$
$$
e_{ii}a(1-e_{ii})xe_{ii}-e_{ii}x(1-e_{ii})ae_{ii}+
e_{ii}d(ii)e_{ii}xe_{ii}-e_{ii}xe_{ii}d(ii)e_{ii}=
$$
$$
e_{ii}a(1-e_{ii})xe_{ii}-e_{ii}x(1-e_{ii})ae_{ii}+0=
$$
$$
e_{ii}a(1-e_{ii})xe_{ii}-e_{ii}x(1-e_{ii})ae_{ii}+e_{ii}ae_{ii}xe_{ii}-e_{ii}xe_{ii}ae_{ii}=
$$
$$
e_{ii}axe_{ii}-e_{ii}xae_{ii}=e_{ii}(ax-xa)e_{ii}
$$
for every $i$.

Hence
$$
\bigtriangleup(x)=ax-xa
$$
for all $x\in SC_{\#}(X,B(H))$. Therefore $\bigtriangleup$ is a
derivation and by \cite{Ol} we may assume that $a\in
SC_{\#}(X,B(H))$. $\triangleright$

\section{The main theorem}

{\bf Theorem 3.1.} {\it Let $M$ be an AW$^*$-algebra of type $I$
and let $\bigtriangleup :M\to M$ be a 2-local derivation. Then
$\bigtriangleup$ is a derivation.}

{\it Proof.} We have that
$$
M=\sum_j^\oplus M_{I_{n_j}},
$$
where $M_{I_{n_j}}$ is an AW$^*$-algebra of type $I_{n_j}$, $n_j$
is a cardinal number for any $j$ and $\sum_j^\oplus M_{I_{n_j}}$
is the C$^*$-sum of the algebras $M_{I_{n_j}}$. Let $x_j\in
M_{I_{n_j}}$ for any $j$ and $x$ be the C$^*$-sum $\sum_j x_j$ of
the elements $x_j$, i.e. $x=\sum_j x_j$. Note that
$\bigtriangleup(x_j)\in M_{I_{n_j}}$ for all $x_j\in M_{I_{n_j}}$.
Hence
$$
\bigtriangleup\vert_{M_{I_{n_j}}}: M_{I_{n_j}}\to M_{I_{n_j}},
$$
$\bigtriangleup$ is a 2-local derivation on $M_{I_{n_j}}$ and by
theorem 2.3 or the theorem in \cite{AA2} $\bigtriangleup$ is a
derivation on $M_{I_{n_j}}$ for $n_j\neq 2$. The case $n_j=2$
follows by the proof of theorem 1 in \cite{AA}.

Let $x$ be an arbitrary element of $M$. Then there exists $d(j)\in
M$ such that $\bigtriangleup(x)=d(j)x-xd(j)$,
$\bigtriangleup(x_j)=d(j)x_j-x_jd(j)$ and
$$
z_j\bigtriangleup(x)=z_j(d(j)x-xd(j))=z_j\sum_i (d(j)x_i-x_id(j))=
$$
$$
d(j)x_j-x_jd(j)=\bigtriangleup(x_j),
$$
for all $j$, where $z_j$ is a unit of $M_{I_{n_j}}$. Hence
$$
\bigtriangleup(x)=\sum_j z_j\bigtriangleup(x)=\sum_j
\bigtriangleup(x_j).
$$
Since $x$ was chosen arbitrarily $\bigtriangleup$ is a derivation
on $M$ by the last equality.

Indeed, let $x,y\in M$. Then
$$
\bigtriangleup(x)+\bigtriangleup(y)=\sum_j
\bigtriangleup(x_j)+\sum_j \bigtriangleup(y_j)= \sum_j
[\bigtriangleup(x_j)+\bigtriangleup(y_j)]=
$$
$$
\sum_j \bigtriangleup(x_j+y_j)=\sum_j z_j\bigtriangleup(x+y)=
\bigtriangleup(x+y).
$$

Similarly,
$$
\bigtriangleup(xy)=\sum_j \bigtriangleup(x_jy_j)=\sum_j
[\bigtriangleup(x_j)y_j+x_j\bigtriangleup(y_j)]=
$$
$$
\sum_j \bigtriangleup(x_j)y_j+\sum_j x_j\bigtriangleup(y_j)=
\sum_j \bigtriangleup(x_j) \sum_j y_j+\sum_j x_j \sum_j
\bigtriangleup(y_j)=
$$
$$
\bigtriangleup(x)y+x\bigtriangleup(y).
$$
Hence $\bigtriangleup$ is a linear operator and a derivation since
$\bigtriangleup$ is homogenous. The proof is complete.
$\triangleright$


\bigskip



\begin{thebibliography}{9999}


\bibitem{S}  \v{S}emrl  P. {\em Local automorphisms and derivations on $B(H)$.} Proc. Amer. Math. Soc. Vol.
125, 2677 - 2680. (1997)

\bibitem{KK}  Kim S.O.,  Kim J.S., {\em Local automorphisms and derivations on $M_n$.} Proc. Amer. Math.
Soc. Vol. 132, 1389 - 1392. (2004)

\bibitem{LW}  Lin Y.,  Wong T. {\em A note on 2-local maps.} Proc. Edinb. Math. Soc. Vol. 49,
701-708. (2006)

\bibitem{AK}  Ayupov Sh.A., Kudaybergenov K.K. {\em 2-local derivations and automorphisms on
$B(H)$.} J. Math. Anal. Appl. Vol. 395, 15-18. (2012)

\bibitem{AA} Ayupov Sh.A., Arzikulov F.N. {\em 2-local derivations on von Neumann algebras of type
I.} arXiv:1112.6236v3 [math.OA] 22 Apr 2014, www.arxiv.org

\bibitem{AKNA} Ayupov Sh.A., Kudaybergenov K.K.,  Nurjanov B.O., Alauatdinov A.K. {\em Local and
2-local derivations on noncommutative Arens algebras.} Mathematica
Slovaca,  64, 423–432. (2014)

\bibitem{AA2} Ayupov Sh.A., Arzikulov F.N., {\em 2-local derivations on semi-finite von Neumann
algebras.} Glasgow Math. Jour. Vol. 56, 9-12. (2014)

\bibitem{AK2} Ayupov Sh.A., Kudaybergenov K.K., {\em 2-local derivations on von Neumann algebras.} Positivity,
DOI 10.1007/s11117-014-0307-3.

\bibitem{Ol} Olesen D., {\em Derivations of AW$^*$-algebras are inner.} Pacific J.
Math. Vol. 53, No. 2, 555-561. (1974)

\bibitem{Bre} Bresar M., {\em Jordan derivations on semiprime rings,} Proc. Amer. Math. Soc. Vol. 104,
1003-1006. (1988)

\bibitem{She} Sherstnev A.N., {\em Methods of bilinear forms in the noncommutative theory of
a measure and an integral.} Moscow, Fizmatlit. (2008, Russian)

\bibitem{Kus} Kusraev A.G., {\em Booleanvalued analysis of involutive
algebras.} Vladikavkaz. (1996, Russian)

\bibitem{Aar} Aarnes J.F. {\em Quasi-states on C*-algebras.} Trans. Am. Math.
Soc. Vol. 149, No 2, 601-625. (1970)


\bibitem{MMS} Matvejchuk M.S. {\em Linearity of a singed measure on a
lattice of ortoprojections.} News of higher educational
institutions. Mathematical series. Vol. 400. No 9, 48-66. (1995,
Russian)


\bibitem{S} Sakai, S.{\em C*-algebras and W*-algebras.} Springer (1971)


\end{thebibliography}
\end{document}